\documentclass[11pt]{amsart}
\usepackage[bottom=1in, right=1in, left=1in, top=1in]{geometry}
\usepackage{graphicx} 
\usepackage{amsmath}

\RequirePackage[OT1]{fontenc}
\RequirePackage{amsthm,amsmath}
\RequirePackage[colorlinks,citecolor=blue,urlcolor=blue]{hyperref}
\RequirePackage{hypernat}
\RequirePackage{comment}
\RequirePackage{graphicx,subfigure,latexsym,amssymb}
\RequirePackage{float,epsfig,multirow,rotating,times}
\RequirePackage{upgreek,wrapfig}
\usepackage{amsthm, amssymb}
\usepackage{tikz}
\usepackage{float}
\usetikzlibrary{matrix}
\usepackage{nicematrix}
\usepackage{arydshln}
\usepackage{blkarray}
\usepackage{xcolor}
\usepackage{amsthm, amssymb}
\newtheorem{theorem}{Theorem}[section]
\newtheorem{proposition}[theorem]{Proposition}

\theoremstyle{definition}
\newtheorem{example}[theorem]{Example}
\theoremstyle{definition}
\newtheorem{definition}[theorem]{Definition}
\theoremstyle{definition}

\theoremstyle{definition}

\usepackage{mathtools} 

\usetikzlibrary{decorations.pathreplacing,angles,quotes}

\renewcommand{\int}{\mathrm{Int}}

\newcommand{\condind}{\perp\!\!\!\perp}

\makeatletter
\renewcommand{\@settitle}{
  \begin{center}
    \normalfont\Large
    \@title
  \end{center}
}
\makeatother

\title[The {GameTheory} package for Macaulay2]{The {GameTheory} package for Macaulay2}

\author[Connelly, Galgano, He, Maletto, Neuhaus, Portakal, Tillmann-Morris, Zhao]{Erin Connelly, Vincenzo Galgano, Zhuang He, Giacomo Maletto, Elke Neuhaus, Irem Portakal, Hannah Tillmann-Morris, Chenyang Zhao}
\date{July 2025}

\AtEndDocument{
    \vspace{5mm}
    \hrule
    \vspace{5mm}
    \noindent {\small (E. Connelly)} \textsc{\footnotesize University of Osnabr\"uck, Osnabr\"uck, Germany} \\ {\small \textit{E-mail address:} {erin.connelly@uni-osnabrueck.de}; ORCID: 0009-0006-0119-3889}\\
    
    \noindent {\small (V. Galgano)} \textsc{\footnotesize Max-Planck Institute of Molecular Cell Biology and Genetics, Dresden, Germany} \\ {\small \textit{E-mail address:} {galgano@mpi-cbg.de}; ORCID: 0000-0001-8778-575X}\\
    
    \noindent {\small (Z. He)} \textsc{\footnotesize University of Turin, Italy} \\ {\small \textit{E-mail address}: { zhuang.he@unito.it }; ORCID: 0000-0003-4082-7355}\\
    
    \noindent {\small (G. Maletto)} \textsc{\footnotesize KTH, Stockholm, Sweden} \\ {\small \textit{E-mail address}: { gmaletto@kth.se }; ORCID: 0000-0003-3211-4689}\\

    \noindent {\small (E. Neuhaus)} \textsc{\footnotesize Max-Planck Institute for Mathematics in the Sciences, Leipzig, Germany} \\ {\small \textit{E-mail address}: { elke.neuhaus@mis.mpg.de }; ORCID: 0009-0000-8478-8891}\\

    \noindent {\small (I. Portakal)} \textsc{\footnotesize Max-Planck Institute for Mathematics in the Sciences, Leipzig, Germany} \\ {\small \textit{E-mail address}: { mail@irem-portakal.de }; ORCID: 0000-0003-3189-5217}\\

    \noindent {\small (H. Tillmann-Morris)} \textsc{\footnotesize Max-Planck Institute for Mathematics in the Sciences, Leipzig, Germany} \\ {\small \textit{E-mail address}: { tillmann@mis.mpg.de }; ORCID: 0000-0001-9541-0296}\\
    
    \noindent {\small (C. Zhao)} \textsc{\footnotesize Imperial College London, United Kingdom} \\ {\small \textit{E-mail address}: { cz2922@ic.ac.uk }; ORCID: 0009-0009-1700-2591}   
}

\begin{document}

\begin{abstract}
    We describe the \texttt{GameTheory} package version~1.0 for computing equilibria in game theory available since version 1.25.05 of \texttt{Macaulay2}. We briefly explain the four equilibrium notions, Nash, correlated, dependency, and conditional independence, and demonstrate their implementation in the package with examples.
 
\end{abstract}
\maketitle
\section{Introduction}\label{sec:Tensors}

We study finite normal-form games and focus on four equilibrium concepts: Nash, correlated, dependency and conditional independence equilibria. Each of these notions possesses rich combinatorial and algebraic structures which can be modeled and explored using \texttt{Macaulay2}. While classical algorithms for computing Nash and correlated equilibria are well established and studied in the literature, such as Lemke-Howson algorithm and linear programming (\cite{nisan2007algorithmic}), our approach uses the capabilities of \texttt{Macaulay2} to provide a unified algebraic and combinatorial framework for all four equilibrium concepts. Therefore, \texttt{GameTheory.m2} enables one to construct and analyze the defining ideals and polyhedral models associated with these equilibria, offering new theoretical insights (\cite{GameTheorySource}). Throughout the package and the paper, indices are \emph{zero-based}: players are numbered as $0,\dots,n-1$ and, if player~$i$ has $d_i$ strategies, they are labeled as $0,\dots,d_i-1$. One caveat is the labeling of vertices in the graphical model, which is \emph{one-based} to ensure compatibility with \verb|GraphicalModels.m2|. This discrepancy does not affect the indexing of variables, and the labels of the vertices of the graphical model can be changed as explained in Section~\ref{sec5}, if necessary. \\
\indent An $n$-player normal-form game is defined by a vector of $n$ tensors $X = (X^{(0)}, \cdots, X^{(n-1)})$ of format $(d_0, \cdots, d_{n-1})$ where $d_i$ is the number of strategies that Player $i$ has. The entry $X^{(i)}_{j_0 \cdots j_{n-1}}$ represents the payoff of Player $i$ when Player $k$ chooses the strategy $j_k \in \{0, \cdots, d_k -1\}$ for $k \in \{0, \cdots, n-1\}$. We call $X^{(i)}$ the payoff tensor (or matrix) of Player $i$ and $X$ is called a $(d_0 \times \cdots \times d_{n-1})$-game. Our package represents an $n$-player game by an ordered list of $n$ tensors of format $(d_0 \times \dots \times d_{n-1})$. We have two ways of constructing a game. The first one is by using \verb|zeroTensor {d_0,...,d_(n-1)}| to create a payoff tensor of the requested format whose entries are all~0, and then modifying those entries to define the normal-form game we have. 
\begin{example}[Bach or Stravinsky]\label{ex: Bach or Stravinsky}
Let us consider the Bach or Stravinsky game in which two players have both two strategies each: going to a Bach concert or a Stravinsky concert. They wish to attend a concert together but have different preferences: the first prefers Bach and the second prefers Stravinsky. Both value being together more than being apart. We define this $2$-player normal-form game $X$ with two $(2 \times 2)$ payoff matrices.
\[
X^{(0)} =
\bordermatrix{
      & \text{B} & \text{S} \cr
\text{B}        & 3 & 0 \cr
\text{S}  & 0 & 2 \cr
}
\qquad \qquad
X^{(1)} =
\bordermatrix{
      & \text{B} & \text{S} \cr
\text{B}        & 2 & 0 \cr
\text{S}  & 0 & 3 \cr
} \ .
\]

\noindent In the package, the game is defined as an ordered list of two tensors as follows.
\vspace{-0.25cm}
{
\footnotesize
\begin{verbatim}

i1 : needsPackage "GameTheory"
i2 : Di = {2,2};
i3 : A = zeroTensor Di;
i4 : B = zeroTensor Di;
i5 : A#{0,0} = 3;  A#{1,1} = 2;
i7 : B#{0,0} = 2;  B#{1,1} = 3;
i9 : X = {A,B};

\end{verbatim}
}

\end{example}

\noindent Alternatively, one may use 
\verb|randomTensor {d_0,...,d_(n-1)}| to create a random tensor of format $(d_0,\cdots,d_{n-1})$ or \verb|randomGame {d_0,...,d_(n-1)}| to return a list of $n$ payoff tensors (a tensor for each player) of format 
$(d_0\times\cdots\times d_{n-1})$ with random entries. For both methods, the ring used by default is $\mathbb{Q}$ and one can change it using the option \verb|CoefficientRing => R|. To see the entries of the tensors, \texttt{peek} is used.

{\footnotesize
\begin{verbatim}
i2 : RG = randomGame {2,2,2}
o2 : {Tensor{...11...}, Tensor{...11...}, Tensor{...11...}}
o2 : List
\end{verbatim}}
\noindent Both constructions above generate tensor objects with the same internal representation as mutable hash tables. A tensor object is equipped with the following properties: \texttt{format}, \texttt{coefficientRing}, and \texttt{indexset}.

\begin{example}\label{exerciseexample}\cite[Exercise 7.2]{mclennan2018advanced}
Let us consider a specific $(2\times2\times2)$-game where three players choose from two strategies $a$ and $b$. In this case, the payoff tensors are generally encoded by two matrices as follows. In each matrix, the first entry in a tuple is the payoff for player one (who picks the row), the second entry is for player two (who picks the column), and the third entry is for player three (who picks which matrix is used). 

\[
\text{Player 3 chooses $a$:}
\quad
\begin{array}{c|cc}
& a & b \\
\hline
a & (1,1,1) & (-5,0,3) \\
b & (0,3,-5) & (0,0,1) \\
\end{array}
\qquad\;
\text{Player 3 chooses $b$:}
\quad
\begin{array}{c|cc}
& a & b \\
\hline
a & (3,-5,0) & (1,0,0) \\
b & (0,1,0) & (0,0,0) \\
\end{array}
\]
\vspace{0.1cm}

{\footnotesize
\begin{verbatim}
i2 : A = zeroTensor {2,2,2};
i3 : A#{0,0,0}=1; A#{0,1,0}=-5; A#{0,0,1}=3; A#{0,1,1}=1;
i7 : B = zeroTensor {2,2,2};
i8 : B#{0,0,0}=1; B#{1,0,0}=3; B#{0,0,1}=-5; B#{1,0,1}=1;
i12 : C = zeroTensor {2,2,2};
i13 : C#{0,0,0}=1; C#{0,1,0}=3; C#{1,0,0}=-5; C#{1,1,0}=1;
i17 : Y = {A,B,C};
\end{verbatim}}

\end{example}
In the following sections, we briefly review the relevant theory and demonstrate how we implement these concepts alongside how the package can be used in practice.

\section{Nash equilibria}
Consider a $(d_0 \times \cdots \times d_{n-1})$-game. Let $p^{(i)}_{j}$ denote the probability that player~$i$ selects pure strategy $j \in \{0,\ldots, d_{n-1}\}$. Then, $\sum_{j=0}^{d_i-1}p^{(i)}_{j}=1$ for each $i$. Therefore a mixed-strategy is a point 
\begin{equation}\label{mixedpoint}
\mathbf{p} =\bigl(p^{(0)}_{0},\dots,p^{(0)}_{d_0-1};\,\dots;\, p^{(n-1)}_{0},\dots,p^{(n-1)}_{d_{n-1}-1}\bigr) 
\end{equation}
that lies in the product of simplexes $\Delta :=\Delta_{d_0-1}\times\cdots\times\Delta_{d_{n-1}-1}$, 
where $\Delta_{d_i-1}$ denotes the probability simplex of dimension $(d_i-1)$.

\begin{definition}
For a $(d_0 \times \cdots \times d_{n-1})$-game, at point~$\mathbf{p}$ in (\ref{mixedpoint}), the \emph{expected payoff} to player~$i$ is given by
\begin{equation*}
\pi_i(\mathbf{p})=\sum_{j_0=0}^{d_0 -1}\cdots\sum_{j_{n-1}=0}^{d_{n-1}}X^{(i)}_{j_0\cdots j_{n-1}}\,p^{(0)}_{j_0}\cdots p^{(n-1)}_{j_{n-1}}.
\end{equation*}
We call $\mathbf{p}$ a \emph{Nash equilibrium} if no player can increase their expected payoff by changing their mixed strategy while assuming the others have fixed mixed strategies.
\end{definition}

\subsection{Totally mixed Nash equilibria and multilinear equations}\label{subsec:polysystem}
A Nash equilibrium $\mathbf{p}$ is called \emph{totally mixed} if $p^{(i)}_{j} > 0$ for all $j \in \{0, \ldots, d_i - 1\}$ and all $i \in \{0, \ldots, n-1\}$, that is, if $\mathbf{p}$ lies in the interior of the product of probability simplices $\Delta$. The set of totally mixed Nash equilibria can be characterized as the set of solutions in the interior of $\Delta$ to the following system of $d_0 + \cdots + d_{n-1} - n$ multilinear equations (see~\cite[Theorem 6.6]{Sturmfels2002}): 
\begin{equation*} \label{multilinear}
\sum_{j_0=0}^{d_0-1} \cdots \widehat{\sum_{j_{i}=0}^{d_{i}-1}} \cdots \sum_{j_{n-1}=0}^{d_{n-1}-1}
\left(
X^{(i)}_{j_0 \cdots j_{i-1} \, k \, j_{i+1} \cdots j_{n-1}}
- X^{(i)}_{j_0 \cdots j_{i-1} \, 0 \, j_{i+1} \cdots j_{n-1}}
\right)
p^{(0)}_{j_0} \cdots p^{(i-1)}_{j_{i-1}} p^{(i+1)}_{j_{i+1}} \cdots p^{(n-1)}_{j_{n-1}} = 0,
\end{equation*}
for all $k \in \{1, \ldots, d_i - 1\}$ and all $i \in \{0, \ldots, n-1\}$.

\begin{example}(Example~\ref{exerciseexample}, continued) 
We first generate the ring $R$ by \texttt{nashEquilibriumRing} which is the polynomial ring over $\mathbb{Q}$ generated by the probability variables $\{p^{(i)}_{j}\}$ for a given game $Y$. Then the method \texttt{nashEquilibriumIdeal} prints the polynomials above with the linear conditions $\sum_{j=0}^{d_i-1}p^{(i)}_{j}=1$. The probability variables $\{p^{(i)}_j\}$ are represented by the symbols \verb|p_{i,j}|.
{\footnotesize
\begin{verbatim}
i18 : R = nashEquilibriumRing Y;
i19 : J = nashEquilibriumIdeal(R, Y);
\end{verbatim}}

\noindent One can verify that this game has a unique totally mixed Nash equilibrium at $((\frac{1}{2},\frac{1}{2}),(\frac{1}{2},\frac{1}{2}),(\frac{1}{2},\frac{1}{2}))$. The equilibrium ideal $J$ has dimension $0$ and degree $2$. The unique totally mixed Nash equilibrium is a double point.
{\footnotesize
\begin{verbatim}
i20 : (dim J, degree J)
o20 = (0, 2)
i21 : decompose J
o21 = {ideal (2p       - 1, 2p       - 1, 2p       - 1, 2p       - 1, 2p      
                {2, 1}        {2, 0}        {1, 1}        {1, 0}        {0, 1}
      --------------------------------------------------------------------------      
      - 1, 2p       - 1)}
             {0, 0}
             
\end{verbatim}}
\noindent We can also perturb the game with a new variable \texttt{e} and define the tensors over \texttt{S = QQ[e]}.
{\footnotesize
\begin{verbatim}
i22 : S = QQ[e];
i23 : A#{0,0,0}=1+e; A#{0,1,0}=-5+e; A#{0,0,1}=3+e; A#{0,1,1}=1+e;
i27 : tensors = {A,B,C};
i28 : R = nashEquilibriumRing tensors;
i29 : R2 = S[gens R];
i30 : J2 = nashEquilibriumIdeal(R2, tensors);
i31 : decompose J2
o31 = {ideal (p       + p       - 1, p       + p       - 1, 6p       - 6p       -
               {2, 0}    {2, 1}       {0, 0}    {0, 1}        {1, 1}     {2, 1}  
      -----------------------------------------------------------------------------                                                                             
      e, 6p       + 6p       + e - 6, (2e - 24)p       + 24p       + e, 4p   
           {1, 0}     {2, 1}                    {0, 1}      {2, 1}        {0,
      -----------------------------------------------------------------------------
                                              2
        p       + 2p       - 6p       + 1, 12p       + (2e - 12)p       + 3)}
      1} {2, 1}     {0, 1}     {2, 1}         {2, 1}             {2, 1}
o31 : List

\end{verbatim}}

\noindent The equilibrium ideal \texttt{J2} corresponds to the game in which the payoff tensor \texttt{A} is perturbed.  
The last generator of the ideal,

\[
12\bigl(p^{(2)}_{1}\bigr)^{2} \;+\; (2e-12)\,p^{(2)}_{1} \;+\; 3,
\]
has a real positive root for \(p^{(2)}_{1}\) if and only if \(e \le 0\). In fact, the game has two totally mixed Nash equilibria when \(-\frac{3}{2}<e< 0\) and a unique totally mixed Nash equilibrium (double point) when \(e=0\). Otherwise, the game has no totally mixed Nash equilibria. 
\end{example}

\subsection{Mixed volume and an upper bound}\label{mixedvol}
By Bernstein's theorem \cite{Bernshtein1975}, the number of isolated complex solutions of multilinear equations from Section~\ref{multilinear} is bounded above by the {\em mixed volume} of the collection of their Newton polytopes. Let $\mathbf{d}=(d_0,\cdots, d_{n-1})$. For a generic game, the $d_i-1$ multilinear polynomials for player $i$ have the same Newton polytope $\Delta^{(i)}$ given by
\[ \Delta^{(i)} := \Delta_{d_{0}-1}\times \Delta_{d_{1}-1} \times \cdots \times \Delta_{d_{i-1}-1} \times \{0\} \times \Delta_{d_{i+1}-1} \times \cdots \times \Delta_{d_{n-1}-1}.\]

\noindent Therefore, the number of isolated totally mixed Nash equilibria of a generic game is bounded by the mixed volume of the list 
\[ \Delta[\mathbf{d}]:=(\Delta^{(0)}, \cdots, \Delta^{(0)},\Delta^{(1)}, \cdots, \Delta^{(1)}, \cdots, \Delta^{(n-1)}, \cdots, \Delta^{(n-1)}),\]
where each $\Delta^{(i)}$ repeats itself $d_i - 1$ times. This $\Delta[\mathbf{d}]$ is implemented by $\texttt{deltaList}$, taking $\mathbf{d}$ as the input.

{\footnotesize
\begin{verbatim}
i2 : DL = deltaList {2,2,2}
o2 = {Polyhedron{...1...}, Polyhedron{...1...}, Polyhedron{...1...}}
o2 : List
i3 : apply(DL, i-> dim i)
o3 = {2, 2, 2}
o3 : List
\end{verbatim}}
\noindent
In \cite[Theorem 3.3]{McKelveyMcLennan1997}, it is shown that, in the case of generic games, this upper bound is achieved for the number of totally mixed Nash equilibria. They establish this by counting block derangements, which is implemented with \texttt{blockDerangements}. The method \texttt{numberTMNE} implements the vector bundle approach of \cite[Thm.\ 2.7]{abo2025vectorbundleapproachnash} for computing the mixed volume of $\Delta[\mathbf{d}]$. Let $\mathbb{P}^{\mathbf{d}} := \prod_{i=0}^{n-1} \mathbb{P}^{d_i-1}$. The vector bundle
\[
E := \bigoplus_{i=0}^{n-1} \mathcal{O}(1,1,\dots, 0^{(i)}, 1,\dots, 1)^{\oplus (d_i-1)}
\]

\noindent
over $\mathbb{P}^{\mathbf{d}}$ has rank $D := \sum_{i=0}^{n-1} (d_i-1)$. The degree of the top Chern class of $E$ coincides with the mixed volume of $\Delta[\mathbf{d}]$, and is equal to the coefficient of the term $\prod_{i=0}^{n-1} h_{i}^{d_i-1}$ in $\prod_{i=0}^{n-1} \left(\sum_{j\neq i} h_j\right)^{d_i-1}$.
The method \texttt{numberTMNE} computes the mixed volume for a given format $\mathbf{d}$ by finding this coefficient. For generic $(2\times 2\times 2)$ games, the mixed volume equals the number of derangements of a $3$-element set, which is $2$.

{\footnotesize
\begin{verbatim}
i4 : numberTMNE {2,2,2}
o4 = 2
\end{verbatim}
}

\noindent
In general, \texttt{numberTMNE} and \texttt{blockDerangements} are the fastest methods, while computing the \texttt{mixedVolume} of \verb|deltaList {3,3,3}| does not terminate within a reasonable amount of time.

{\footnotesize
\begin{verbatim}
i5 : elapsedTime numberTMNE {3,3,3}
 -- .0194357s elapsed
o5 = 10
i6 : elapsedTime #(blockDerangements {3,3,3})
 -- .0481826s elapsed
o6 = 10
\end{verbatim}
}
\section{Correlated equilibria}\label{correlatedequilibria}
The concept of correlated equilibria was introduced by Robert Aumann in 1974 \cite{aumann1974subjectivity}. Unlike Nash equilibria, where players choose strategies independently, Aumann's concept allows for strategic coordination through shared signals, expanding the possibilities for improved outcomes in normal-form games. We consider joint probability distributions instead of disjoint mixed strategies for each player. The variable $p_{j_0 j_1 \cdots j_{n-1}}$ represents the probability the $i$-th player chooses the strategy $j_i \in \{0,\ldots, d_i-1\}$, for all $i\in \{0,\cdots,n-1\}$. 
To represent probabilities, the probability tensor needs to live in the probability simplex $\Delta_{(d_0 \cdots d_{n-1})-1}$.
Such variables define the ring of joint probability distributions of a game, which is implemented by \texttt{probabilityRing}: it takes as input the format of the game $(d_0,\ldots, d_{n-1})$, and returns as output a polynomial ring with a tensor of variables $p_{j_0 \cdots j_{n-1}}$. For instance, for a game format $(2,1, 2)$ one obtains:

{\footnotesize
\begin{verbatim}
i2 : Di = {2,1,2};
i3 : PR = probabilityRing Di
o3 = PR
o3 : PolynomialRing
i4 : pairs PR#"probabilityVariable"
o4 = {(indexes, {{0, 0, 0}, {0, 0, 1}, {1, 0, 0}, {1, 0, 1}}), (format, {2, 1, 2}),
      -------------------------------------------------------------------------------
      ({0, 0, 0}, p         ), ({1, 0, 0}, p         ), ({0, 0, 1}, p         ), 
                   {0, 0, 0}                {1, 0, 0}                {0, 0, 1} 
      -------------------------------------------------------------------------------
      ({1, 0, 1}, p         ), (coefficients, PR)}
                   {1, 0, 1}      
o4 : List
\end{verbatim}
}

\noindent It is important to note that for $2$-player games, this output is not the same as the output for probabilities in the \texttt{nashEquilibriumRing}, although these print the same.
The base field of the probability ring is set to $\mathbb Q$, and the tensor of probability variables to $p$. It is possible to change the base field and the name of the tensor by using the optional arguments \texttt{CoefficientRing} and \texttt{ProbabilityVariableName}. 

{\footnotesize
\begin{verbatim}
i5 : PR2 = probabilityRing (Di, CoefficientRing=>RR, ProbabilityVariableName=>"q");
i6 : coefficientRing PR2;
i7 : pairs PR2#"probabilityVariable";
\end{verbatim}
}

\noindent While Nash equilibria are generically finite in number, the set of correlated equilibria is typically infinite. In fact, the correlated equilibria form a convex polytope within the probability simplex $\Delta_{(d_0 \cdots d_{n-1})-1}$.
\begin{definition}[\cite{aumann1987correlated},\cite{brandenburg2025combinatorics}]
   A point $p = (p_{j_0 \cdots j_{n-1}}) \in \Delta_{(d_0 \cdots d_{n-1}) -1}$ is a \emph{correlated equilibrium} if and only if 
	\begin{equation}\label{eq: correlated equilibria}
    \sum_{j_0=0}^{d_0-1} \cdots \widehat{\sum_{j_{i}=0}^{d_{i}-1}} \cdots \sum_{j_{n}=0}^{d_{n}-1}
    \left( X^{(i)}_{j_0\cdots j_{i-1} k j_{i+1} \cdots j_{n-1}} - X^{(i)}_{j_0 \cdots j_{i-1} l j_{i+1} \cdots j_{n-1}} \right)
    p_{j_0 \cdots j_{i-1} k j_{i+1} \cdots j_{n-1}} \geq 0.
\end{equation}
for all $k, l \in \{0, \ldots, d_i-1\}$, and for all $i \in \{0, \ldots, n-1\}$. The linear inequalities in \eqref{eq: correlated equilibria} together with the linear constraints
\[
    p_{j_0 \dots j_{n-1}} \geq 0 \text{ for } j_i \in \{0, \ldots, d_i-1\},\ i \in \{0, \ldots, n-1\}, \quad \text{ and }\ \ 
    \sum_{j_0=0}^{d_0-1} \cdots \sum_{j_{n-1}=0}^{d_{n-1}-1} p_{j_0 \dots j_{n-1}} = 1
\]
define the \emph{correlated equilibrium polytope} $P_X$ of the game $X = (X^{(0)},\dots,X^{(n-1)})$.

\end{definition}

\begin{example}[Bach or Stravinsky, Example~\ref{ex: Bach or Stravinsky} continued] The Nash equilibria of this game consist of two pure Nash equilibria where both choose to go to Bach concert $(1,0,0,0)$ and to go to Stravinsky concert $(0,0,0,1)$. One also has a unique totally mixed Nash equilibrium $\left( \frac{6}{25},\; \frac{9}{25},\; \frac{4}{25},\; \frac{6}{25} \right)$. Now, one can consider the situation where both players agree to use an external randomizing device (such as a coin flip) to coordinate their choices: choosing (Bach, Bach) if heads and (Stravinsky, Stravinsky) if tails. The coin flip mediates this coordination, allowing them to randomize fairly over the two preferred outcomes, which are Nash equilibria, resulting in a correlated equilibrium that is Pareto optimal, that neither player has an incentive to deviate from.
\end{example}

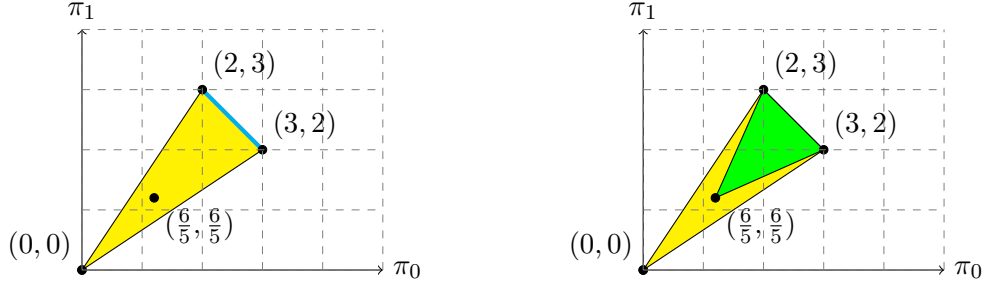
\begin{figure}[H]
\begin{tikzpicture}[scale=0.8]
    \draw[->] (0,0) -- (5,0) node[right] {$\pi_0$};
    \draw[->] (0,0) -- (0,4) node[above] {$\pi_1$};
    
    \filldraw[fill=yellow, draw=black] (0,0) -- (2,3) -- (3,2) -- cycle;

    \filldraw[fill=cyan, draw=cyan, ultra thick] (2,3) -- (3,2) -- cycle;
    
    \filldraw (3,2) circle (2pt) node[above right] {$(3,2)$};
    \filldraw (2,3) circle (2pt) node[above right] {$(2,3)$};
    \filldraw (6/5,6/5) circle (2pt) node[below right] {$\color{black}(\frac{6}{5},\frac{6}{5})$};
    \filldraw (0,0) circle (2pt) node[above left] {$(0,0)$};
    \filldraw (0,0) circle (2pt) node[above left] {};
    
    \draw[help lines, dashed, gray] (0,0) grid (5,4);
\end{tikzpicture}
\hspace{1.5cm}
\begin{tikzpicture}[scale=0.8]
    \draw[->] (0,0) -- (5,0) node[right] {$\pi_0$};
    \draw[->] (0,0) -- (0,4) node[above] {$\pi_1$};
    
    \filldraw[fill=yellow, draw=black] (0,0) -- (2,3) -- (3,2) -- cycle;

    \filldraw[fill=green, draw=black] (2,3) -- (3,2)  -- (6/5,6/5) -- cycle;
    
    \filldraw (3,2) circle (2pt) node[above right] {$(3,2)$};
    \filldraw (2,3) circle (2pt) node[above right] {$(2,3)$};
    \filldraw (6/5,6/5) circle (2pt) node[below right] {$\color{black}(\frac{6}{5},\frac{6}{5})$};
    \filldraw (0,0) circle (2pt) node[above left] {$(0,0)$};
    \filldraw (0,0) circle (2pt) node[above left] {};
    
    \draw[help lines, dashed, gray] (0,0) grid (5,4);
\end{tikzpicture}
\caption{\small The three points are the payoffs of the players for Nash equilibria. The yellow triangle is the image of $\Delta_3$ in the payoff space. The blue line corresponds to the randomization between two pure Nash equilibria. The green triangle is the image of the correlated equilibrium polytope of the game.}
\end{figure}
\noindent The package uses \texttt{Polyhedra.m2} for studying the polyhedral geometry of the correlated equilibrium polytope. In this case, the correlated equilibrium polytope $P_X$ is a double pyramid over a triangle that is the convex hull of three Nash equilibria of $X$.
\vspace{0.2cm}
\begin{center}
\begin{minipage}[t]{0.48\textwidth}
\footnotesize
\begin{verbatim}
i10 : CE = correlatedEquilibria X
o10 : Polyhedron
i11 : dim CE
o11 = 3
i12 : vertices CE
o12 = | 1 0 2/7 3/8 6/25 |
      | 0 0 3/7 0   9/25 |
      | 0 0 0   1/4 4/25 |
      | 0 1 2/7 3/8 6/25 |
      
               4       5
o12 : Matrix QQ  <-- QQ
\end{verbatim}
\end{minipage}%
\hfill
\begin{minipage}[t]{0.48\textwidth}
\footnotesize
\begin{verbatim}
i13 : facets CE
o13 = (| 0  -1 0  0  |, 0)
       | -3 2  0  0  |
       | 0  0  -1 0  |
       | -2 0  3  0  |
       | 0  2  0  -3 |
       | 0  0  3  -2 |
o13 : Sequence      
i14 : fVector CE  
o14 = {5, 9, 6, 1}
o14 : List
\end{verbatim}
\end{minipage}
\end{center}
\section{Dependency equilibria}
While a correlated equilibrium ensures that no player can increase their payoff by deviating from the recommended joint probability distribution, in a dependency equilibrium each player's conditional expected payoff does not depend on their choice of strategy. Dependency equilibria were introduced by Spohn \cite{spohn2007dependency}, and their geometry was studied by Portakal-Sturmfels \cite{portakal2022geometry} who introduced the related notion of Spohn varieties.

\begin{definition}\cite{portakal2022geometry, spohn2007dependency} The \emph{conditional expected payoff} of player $i$ conditioned on them choosing pure strategy $k$ is given by

    $$
    \mathbb{E}^{(i)}_k (p) := 
     \sum\limits_{j_0=0}^{d_0-1} \cdots \widehat{\sum\limits_{j_i=0}^{d_i-1}} \cdots \sum\limits_{j_{n-1}=0}^{d_{n-1}-1} X^{(i)}_{j_0\cdots k \cdots j_{n-1}} \frac{p_{j_0\cdots k \cdots j_{n-1}}}{p_{+\ldots+k+ \ldots +}},
    $$
where $p_{+\ldots+k+\ldots+}$ is the overall probability of player $i$ choosing strategy $k$. For a joint probability distribution $p \in \Delta_{d_0 \cdots d_{n-1}-1}^{\circ}$, for which the conditional expected payoffs are well-defined, we call $p$ a \emph{dependency equilibrium} if
    $$
     \mathbb E^{(i)}_k (p) \geq \mathbb E^{(i)}_{k'} (p)   
    $$
    for all players $i \in \{ 0, \ldots n-1 \}$ and all pure strategies $k,k' \in \{ 0, \ldots, d_i-1 \}$ of player $i$.
    In particular, 
    $
     \mathbb E^{(i)}_k (p) = \mathbb E^{(i)}_{k'} (p).
    $

\noindent These probability distributions then lie on the \emph{Spohn variety} $\mathcal V_X$, cut out by the vanishing of the $2 \times 2$-minors of the \emph{Spohn matrices}
    $$M_i(p) := 
    \begin{bmatrix}
    \vdots & \vdots \\
    p_{+\ldots+k+ \ldots +}& \sum\limits_{j_0=0}^{d_0-1} \cdots \widehat{\sum\limits_{j_i=0}^{d_i-1}} \cdots \sum\limits_{j_{n-1}=0}^{d_{n-1}-1} X^{(i)}_{j_0\cdots k \cdots j_{n-1}} p_{j_0\cdots k \cdots j_{n-1}} \\
    \vdots & \vdots
    \end{bmatrix}
    \in \mathbb{R}^{d_i \times 2}.
    $$
\end{definition}

\begin{theorem}[Theorem 6, \cite{portakal2022geometry}]
If the payoff tables $X$ are generic, then the Spohn variety $\mathcal{V}_X$ is irreducible of
codimension $d_0 + d_1 + \ldots + d_{n-1} - n$ and degree $d_0d_1 \cdots d_{n-1}$. The intersection of $\mathcal{V}_X$ with the Segre variety in the open simplex $\Delta_{d_0 \cdots d_{n-1}-1}^{\circ}$ is precisely the set of totally mixed Nash equilibria for $X$.
\end{theorem}

Similarly to correlated equilibria, the underlying space of joint probability distributions has coordinate ring given by a \texttt{probabilityRing}.
The (list of) Spohn matrices and the ideal defining the Spohn variety can be computed using the commands \texttt{spohnMatrices} and \texttt{spohnIdeal}. These take as input the underlying \texttt{probabilityRing} as well as a list of payoff tensors. After fixing the format of the game, we can compute \texttt{spohnMatrices}
and \texttt{spohnIdeal} for a randomly generated game of that format.

{\footnotesize
\begin{verbatim}
i2 : Di = {2,2,2};
i3 : PR = probabilityRing Di;
i4 : X = randomGame Di;

i5 : spohnM = spohnMatrices(PR,X);
i6 : spohnI = spohnIdeal(PR,X);

i7 : ( codim spohnI , codim spohnI == sum Di - 3 )
o7 = (3, true)
o7 : Sequence

i8 : ( degree spohnI , degree spohnI == product Di )
o8 = (8, true)
o8 : Sequence

i9 : isPrime spohnI
o9 = true
\end{verbatim}
}

\noindent The Spohn variety $\mathcal{V}_X$ can alternatively be defined as the locus of joint probability distributions $p$ such that $K_X(z)\cdot p=0$ for some $z\in(\mathbb{P}^1)^d$, where $K_X(z)$ is the \textit{Konstanz matrix} defined by
\begin{equation}\label{eq: def Konstanz}
K_X(z)=\begin{pmatrix}K_{X,0}(z_0)\\\vdots\\ K_{X,n-1}(z_{n-1})\end{pmatrix},\qquad K_{X,i}(z_i)\cdot p=M_i(p)\cdot z_i^T\quad\text{for $i=0,\dots,n-1$},
\end{equation}
where the tensor $p$ is vectorized as a column. 
When working over the affine chart $z_i=(k_i:1)$, one usually also writes $K_X(k)$.
In the case of generic payoff tensors, the Konstanz matrix is also used in the construction of an inverse to the birational morphism $\mathcal{V}_X\dashrightarrow\mathbb{P}^{D-1}\times(\mathbb{P}^1)^n$ with $D=d_0\cdots d_{n-1}-(d_0+\dots+d_{n-1})$, as described in the proof of the following result.

\begin{theorem}[Theorem 8, \cite{portakal2022geometry}]
If $n = d_0 = d_1 = 2$, then the Spohn variety $\mathcal{V}_X$ is an elliptic curve. In all other
cases, the Spohn variety $\mathcal{V}_X$ is rational, represented by a map onto $(\mathbb{P}^1)^n$ with linear fibers.
\end{theorem}

\noindent  The Konstanz matrix $K_X(k)$ can be computed using the command \texttt{konstanzMatrix}, taking as input the base \texttt{probabilityRing} and a list of payoff tensors.
Optionally, one can adapt the name of the new variable $k$ via \texttt{KonstanzVariableName => k}.

\begin{example}[Bach or Stravinsky, Example \ref{ex: Bach or Stravinsky} continued]
Here, we compute \texttt{spohnMatrices}, \texttt{spohnIdeal}, and \texttt{konstanzMatrix} for this game. We also show that the corresponding Spohn variety is reducible and exhibit a primary decomposition of the ideal. One can compare the outputs below with \cite[Example 5]{portakal2022geometry}.
{\footnotesize
\begin{verbatim}
i15 : PR = probabilityRing Di;
i16 : spohnM = spohnMatrices(PR,X)
o16 = {| p_{0, 0}+p_{0, 1} 3p_{0, 0} |, | p_{0, 0}+p_{1, 0} 2p_{0, 0} |}
       | p_{1, 0}+p_{1, 1} 2p_{1, 1} |  | p_{0, 1}+p_{1, 1} 3p_{1, 1} |
o16 : List

i17 : spohnI = spohnIdeal(PR,X)
o17 = ideal (- 3p      p       - p      p       + 2p      p      , 
                 {0, 0} {1, 0}    {0, 0} {1, 1}     {0, 1} {1, 1}
      -------------------------------------------------------------
        - 2p      p       + p      p       + 3p      p      )
            {0, 0} {0, 1}    {0, 0} {1, 1}     {1, 0} {1, 1}

o17 : Ideal of PR 

i18 : primaryDecomposition spohnI
o18 = {ideal (p      , p      ), ideal (2p       - 3p       - p      , p       - p      ), 
               {1, 1}   {0, 0}            {0, 1}     {1, 0}    {1, 1}   {0, 0}    {1, 1}   
      ------------------------------------------------------------------------------------
      ideal (2p       + 3p      , 3p     p       + p      p       + 3p      p      )}
               {0, 1}     {1, 0}    {0,0} {1, 0}    {0, 0} {1, 1}     {1, 0} {1, 1}      
o18 : List

i19 : K = konstanzMatrix(PR,X)
o19 = | k_0-3 k_0 0   0     |
      | 0     0   k_0 k_0-2 |
      | k_1-2 0   k_1 0     |
      | 0     k_1 0   k_1-3 |
                         4                 4
o19 : Matrix (PR[k ..k ])  <-- (PR[k ..k ])
                  0   1             0   1    
\end{verbatim}
}

\noindent We now check the definition of the Konstanz matrix \eqref{eq: def Konstanz} for the Bach or Stravinsky game. We first need to vectorise the tensor $p$ of probability variables.

{\footnotesize
\begin{verbatim}
i20 : p = PR#"probabilityVariable";
i21 : J = enumerateTensorIndices Di;
i22 : P = vector(apply(J, j -> p#j));
i23 : K*P
o23 = | (p_{0, 0}+p_{0, 1})k_0-3p_{0, 0} |
      | (p_{1, 0}+p_{1, 1})k_0-2p_{1, 1} |
      | (p_{0, 0}+p_{1, 0})k_1-2p_{0, 0} |
      | (p_{0, 1}+p_{1, 1})k_1-3p_{1, 1} |

                  4
o23 : (PR[k ..k ])
           0   1

i24 : RHS = (spohnM_0 * vector{k_0, -1}) || (spohnM_1 * vector{k_1, -1});

                  4
o24 = (PR[k ..k ])
           0   1
i25 : flatten entries (K*P) == flatten entries RHS
o25 = true
\end{verbatim}
}

\end{example}

\begin{example}

As pointed out above, the Spohn variety of the Bach and Stravinsky game is reducible, and thus not an elliptic curve. However, if we slightly perturb the game (e.g. with \verb|A#{0,1}=1|), we gain back genericity (cf. \cite[Example 4.12]{kidambi2025elliptic}).

{\footnotesize
\begin{verbatim}
i2 : Di = {2,2};
i3 : PR = probabilityRing Di;
i4 : A = zeroTensor Di; 
i5 : A#{0,0}=3; A#{0,1}=1; A#{1,1}=2;
i8 : B = zeroTensor(Di);
i9 : B#{0,0}=2; B#{1,1}=3;
i11 : X = {A,B};

i12 : spohnI = spohnIdeal(PR,X)
o12 = ideal (- 3p      p       - p      p       - p      p       + p      p      ,
                 {0, 0} {1, 0}    {0, 1} {1, 0}    {0, 0} {1, 1}    {0, 1} {1, 1}  
    -------------------------------------------------------------------------------
     - 2p      p       + p      p       + 3p      p      )
         {0, 0} {0, 1}    {0, 0} {1, 1}     {1, 0} {1, 1}
o12 : Ideal of PR
i13 : isPrime spohnI
o13 = true
\end{verbatim}
}

\end{example}

\section{Conditional independence equilibria}\label{sec5}
In this section, we introduce conditional independence statements to the graphical model of a game, and present the methods available in this package to compute the Spohn conditional independence (CI) variety. The Spohn CI varieties were introduced in \cite[Section 6]{portakal2022geometry} and studied further in \cite{portakal2024nash,portakal2025game}. Note that the vertices of the graphical model are numbered one-based to ensure compatibility with \verb|GraphicalModels.m2|. In particular, this does not affect the indices of the generators of \verb|probabilityRing|. The labeling of the vertices can be adjusted, as demonstrated in Example~\ref{ex: relabeling}, when required.

\subsection{Conditional independence models}

Given a game $X$ with $n$ players, one may wish to model conditional dependencies between the players' strategies. 
If $A$, $B$ and $C$ are disjoint subsets of the $n$ players, we write $X_A\condind X_B\, |\, X_C$ if the strategies of players $A$ are independent of the strategies of the players $B$ given the strategies of $C$.

\begin{proposition}
(\cite{Sul18}) For a subset $A\subset \{0, \cdots, n-1\}$ let $R_A:=\Pi_{a\in A} [d_a]$ be the set of pure strategies for $A$. Then the conditional independence $X_A\condind X_B|X_C$ holds if and only if
\begin{equation}\label{eq:cond_ind}
    p_{i_Ai_Bi_C}+p_{j_Aj_Bi_C}-p_{i_Aj_Bi_C}+p_{j_Ai_Bi_C}=0
\end{equation}
for all $i_A,j_A\in R_A,~i_B,j_B\in R_B,~i_C,j_C\in R_C$.
\end{proposition}

\noindent Given a set $\mathcal C$ of conditional independence statements, we define the discrete conditional independence model $\mathcal M_\mathcal C$ as the variety cut out by all equations of the form \eqref{eq:cond_ind}. 
The method \texttt{ciIdeal} produces the conditional independence ideal as an ideal of the \texttt{probabilityRing}, given a list of conditional independence statements, where each conditional independence statement $X_A\condind X_B\, |\, X_C$ is encoded as a list \texttt{\{A,B,C\}}.

\begin{example}\label{ex: line graph cimodel}
Suppose that $X$ is a $(2 \times 2 \times 2)$-game where the strategies of Players 1 and 3 are independent given Player 2's strategy. 

{\footnotesize
\begin{verbatim}
i2 : FF = ZZ/32003;
i3 : Di = {2,2,2};
i4 : PR = probabilityRing (Di, CoefficientRing => FF);
i5 : Stmts = {{{1}, {3}, {2}}};

i6 : I = ciIdeal (PR, Stmts)
o6 = ideal (- p         p          + p         p         , - p         p
              {0, 0, 1} {1, 0, 0}    {0, 0, 0} {1, 0, 1}     {0, 1, 1} {1, 1, 0}
      --------------------------------------------------------------------------
      + p         p          )
          {0, 1, 0} {1, 1, 1}
o6 : Ideal of PR
\end{verbatim}
}
\end{example}

\noindent The method \texttt{ciIdeal} works by calling the \texttt{conditionalIndependenceIdeal} method from the \texttt{GraphicalModels} package and then constructing an isomorphism to map the ideal from the \texttt{markovRing} of the \texttt{GraphicalModels} package to the \texttt{probabilityRing} of our package. 
To allow further utility in interacting with the \texttt{GraphicalModels} package, we have made these methods callable by the user. The method \texttt{toMarkovRing} takes a ring created with \texttt{probabilityRing} and outputs the canonically isomorphic \texttt{markovRing} from the \texttt{GraphicalModels} package; and both the methods \texttt{mapToMarkovRing} and \texttt{mapToProbabilityRing} produce the isomorphism between the two rings as a \texttt{RingMap}.
The default variable names of the ring produced by \texttt{toMarkovRing} are $p$ unless the variable names of the input ring are $p$, in which case the variable names of the output ring become $q$.

\begin{example}
Continuing with the ring from Example~\ref{ex: line graph cimodel}:

{\footnotesize
\begin{verbatim}
i7 : toMarkovRing(PR)
o7 = FF[q     ..q     ]
         1,1,1   2,2,2
o7 : PolynomialRing
\end{verbatim}
}
\end{example}

\subsection{Discrete undirected graphical models and Markov property}

Given a game with $n$ players, one can encode the conditional independence of their strategies via a graph with $n$ vertices. Each vertex is associated to a player, and the conditional independence statements are extracted via its \textit{Markov property}. More precisely, a pair of vertices $(a,b)$ on a graph $G = ([n], E)$ is said to be \textit{separated by} a subset $C \subseteq [n] \setminus \{a,b\}$ of vertices if every path from $a$ to $b$ contains a vertex $c \in C$. 
If $A, B, C \subseteq [n]$ are disjoint subsets of vertices, we say that $C$ separates $A$ and $B$ if $a$ and $b$ are separated by $C$ for all $a \in A$ and $b \in B$.
The \textit{global Markov property} associated to $G$ is the set $\text{global}(G)$ of CI statements of the form $X_A\condind X_B\, |\, X_C$ for all triples $(A, B, C)$ of disjoint subsets of vertices such that $C$ separates $A$ and $B$. 

There are also pairwise and local Markov properties associated to a graph $G$ -- the \texttt{GraphicalModels} package has methods to compute all of them. 
However, the global Markov property is integrated into the \texttt{ciIdeal} method -- the input can be given as an undirected graph rather than a list of CI statements.

\begin{example}
The line graph connecting three vertices has global Markov property consisting of the single CI statement $X_{1}\condind X_{3}\, |\, X_{2}$. 
Continuing with Example~\ref{ex: line graph cimodel}:
{\footnotesize
\begin{verbatim}
i7 : G = graph ({{1,2},{2,3}});
i8 : I == ciIdeal (PR, G)
o8 = true
\end{verbatim}
}
\end{example}

\noindent The user also has the option to rename the vertices of the graph/players of the game.

\begin{example}\label{ex: relabeling} Let the respective names of the players from Example~\ref{ex: line graph cimodel} be Alice, Bob, and Claire. Then:

{\footnotesize
\begin{verbatim}
i7 : PlayerNames = {Alice, Bob, Claire};
i8 : G = graph ({{Alice, Bob}, {Bob, Claire}});
i9 : I == ciIdeal (PR, G, PlayerNames)
o9 = true
\end{verbatim}
}
\noindent One can also relabel the vertices of the graph as follows:
{\footnotesize
\begin{verbatim}
i10 : LabelNames = {0,1,2};
i11 : G = graph ({{0,1},{1,2}});
i12 : I == ciIdeal (PR, G, LabelNames)
o12 = true
\end{verbatim}
}
\end{example}

\subsection{The Spohn CI variety}
Finally, we incorporate the conditional independence conditions into our game~$X$. 
Let $W$ be the union of all hyperplanes of the form $\{p_{j_0\ldots,j_{n-1}}=0\}$ and $\{p_{+\ldots +}=0\}$. Let $\mathcal{M}_\mathcal{C}^{\text{par}}$ be the parametrized discrete undirected graphical model
associated to $G$ i.e.\ $\mathcal{M}_{\mathcal C}$ with $W$ removed (see \cite[Section 2.2]{portakal2025game}). We then define the \textit{Spohn CI variety} as
\begin{equation*}
    \mathcal{V}_{X,\mathcal C}=\overline{(\mathcal{V}_X\cap \mathcal{M}_\mathcal{C}^{\text{par}})\backslash W}.
\end{equation*}
The method \texttt{spohnCI} produces the ideal defining the Spohn CI variety $\mathcal{V}_{X,\mathcal C}$ given the inputs of a game $X$ and either a list of conditional independence statements \texttt{Stmts} or a graph $G$.

\begin{example}
For example, let $X$ be a random $(2\times 2\times 2)$-game where the strategies of Players 1 and 3 are independent given Player 2's strategy.
{\footnotesize
\begin{verbatim}
i7 : X = randomGame(Di, CoefficientRing => FF);
i8 : SpohnCI = spohnCI (PR, X, Stmts);
o8 : Ideal of PR
\end{verbatim}
}
\end{example}

\noindent The method \texttt{spohnCI} works by producing the ideal of the Spohn variety $\mathcal{V}_X$ by calling \texttt{spohnIdeal}, then applying a separately-callable method \texttt{intersectWithCImodel} to this ideal to produce the ideal defining $\mathcal{V}_{X,\mathcal C}$. The method \texttt{intersectWithCImodel} takes an ideal $V$ of the \texttt{probabilityRing} and a conditional independence model (in the form of a list of CI statements \texttt{Stmts} or a graph \texttt{G}) and outputs the saturation of the sum of $V$ and the conditional independence ideal $I$ by the hyperplanes $\{p_{j_0\ldots,j_{n-1}}=0\}$ and $\{p_{+\ldots +}=0\}$.

\begin{example} 
The \texttt{spohnCI} method can be broken down as follows:
{\footnotesize
\begin{verbatim}
i9 : V = spohnIdeal(PR, X);
o9 : Ideal of PR
i10 : SpohnCI == intersectWithCImodel (V, Stmts)
o10 = true 
\end{verbatim}
}  
\end{example}

\noindent If the intersection between the conditional independence variety and the Spohn variety is empty then the method terminates early. Otherwise, we proceed to saturate the intersection by the hyperplanes we wish to remove. This saturation takes a significant amount of time, but we have found experimentally that this computation runs more efficiently if we saturate our two ideals before intersecting, and then repeat this saturation afterwards. The former saturation runs significantly faster, so we recommend the use of the \texttt{Verbose} option to track the saturation process in the methods \texttt{spohnCI} and \texttt{intersectWithCImodel} when running lengthy computations; this can serve as an estimate for the length of time remaining. The method uses the Bayer strategy for each saturation step.

\begin{example}
The \texttt{Verbose} option applied in the case of our example results in an output that looks as follows.
{\footnotesize
\begin{verbatim}
i11 : intersectWithCImodel(V, Stmts, Verbose => true);
Completed step 1 of saturating CI ideal
Completed step 1 of saturating input ideal
Completed step 2 of saturating CI ideal
\end{verbatim}
$\vdots$
\begin{verbatim}
Completed step 9 of saturating input ideal
Completed step 1 of saturating sum
\end{verbatim}
$\vdots$
\begin{verbatim}
Completed step 8 of saturating sum

o11 : Ideal of PR
\end{verbatim}
}  
\end{example}

\subsection*{Acknowledgement}
We thank Ben Hollering and Mahrud Sayrafi for their support during the workshop \textquotedblleft Macaulay2 in the Sciences\textquotedblright \space at MPI-MiS Leipzig, where the development of this package began. We are grateful to Lars Kastner and Luca Sodomaco for their contribution to this package. VG is member of the italian “National Group for Algebraic and Geometric Structures, and their Applications” (GNSAGA-INdAM). Z. He has been funded by the European Union - Next Generation EU, Mission 4, Component 1, CUP D53D23005860006, within the project PRIN 2022L34E7W “Moduli spaces and birational geometry”.

\bibliographystyle{plain}
\bibliography{references}

\end{document}